# Approximation of Some Nonlinear Fractional Order BVPs by Weighted Residual Methods

Umme Ruman[1] and Md. Shafiqul Islam[2,*]

[1]Department of Computer Science & Engineering, Green University of Bangladesh, Dhaka, Bangladesh
[2]Department of Applied Mathematics, University of Dhaka, Dhaka, Bangladesh

*Corresponding author: mdshafiqul@du.ac.bd

## Abstract

To extract the approximate solutions in the case of nonlinear fractional order differential equations with the homogeneous and nonhomogeneous boundary conditions, the weighted residual method is embedded here. We exploit three methods such as Galerkin, Least Square, and Collocation for the efficient numerical solution of nonlinear two-point boundary value problems. Some nonlinear cases are examined for observing the maximum absolute errors by the considered methods, demonstrating the accuracy and reliability of the present technique using the modified Legendre and modified Bernoulli polynomials as weight functions. The mathematical formulations and computational algorithms are more straightforward and uncomplicated to understand. Absolute errors and the graphical representation reflect that our method is more accurate and reliable.

## 1. Introduction

The differential equation of fractional order is acquainted with several areas of biomedical science, fluid dynamics, and engineering fields. Arising in these fields, fractional differential equations (FDE) with boundary conditions of linear or nonlinear problems have been attempted for solving analytically [1] or numerically [2, 3]. The generalized properties of fractional derivatives and integrals and their applications were established by Delkhosh [4]. Linear fractional order two-point boundary value problems have been solved by many methods such as Adomian decomposition [5], Sinc-Galerkin [6], Cubic spline solution [7]. In [8], Cubic B-Spline wavelet collocation method has been established while Legendre wavelet approximation method has been developed in the literature [9]. Further, Collocation-shooting method and Galerkin WRM have been depicted in the works [10] and [11], respectively.

In the past few years, most of the works have been devoted to solving the nonlinear initial and boundary value problems for fractional order differential equations by various methods, such as, Elzaki and Chamekh introduced a New Decomposition Method [12], the variational iteration method was applied by Momani [13], Runge-Heun-Kutta method by Harker [14],

Sinc- Galerkin method by El-Gamel [15], homotopy analysis method by El-Ajou et al [16], Legendre wavelet operational matrix method by Secer et al [17], Pseudospectral method by Delpasand [18], the operational matrix method by Lakestani et al [19], The Least-square method by Fadhel and Jameel [20]. In [21], Bai and Lu investigated the existence and multiplicity of the positive solutions for nonlinear BVP of fractional order.

Recently, Schauder fixed point theorem was used to investigate the existence of solutions to the initial value problem for nonlinear fractional differential equations involving Caputo sequential fractional derivative by Ye and Huang [22].

From the above literature review on numerical solutions of nonlinear fractional differential equations, we are motivated to find an efficient and reliable method for fractional BVP. Therefore, the objective of this paper is to introduce the application of the weighted residual method for finding the approximation to the two-point BVP of fractional order which is as follows:

$$D_*^\alpha x(t) = f\big(t, x(t), x'(t), x''(t), \ldots, x^n(t)\big), \qquad a \le t \le b$$

subject to the boundary conditions:

$$x^{(i)}(a) = a_i, \qquad y^{(i)}(b) = b_i, \quad i = \{0, 1, 2, \ldots, k-1\}$$

where $D_*^\alpha x(t)$ is the fractional derivative of order $\alpha$ of $x(t)$ in the Caputo sense and $k \epsilon N$.

The goal of the proposed research work is to use the Galerkin, Least Square, and Collocation Weighted Residual Methods for solving nonlinear fractional order boundary value problems.

In order to prepare this research work we organise as follows. Some basic definitions of fractional derivatives and notations of fractional calculus are defined in Section 2. The mathematical formulation of three proposed methods for nonlinear fractional order differential equations are given elaborately in Section 3. The numerical solutions to the specific problems and the comparison of the absolute errors of different methods are displayed in tabular form in section 4, and finally the conclusion and references are appended.

## 2. Preliminaries

Fractional Derivatives: In the study of derivatives, we need to solve the non-integer derivatives which contained in the fractional differential equation. There are several definitions of fractional derivatives such as Grunwald-Letnikov fractional derivatives, Riemann-Liouville fractional derivative and Caputo fractional derivative. The effectiveness of these definitions has been utilized but there were some limitations that deal with the initial conditions. To overcome this deficiency, we have to establish the solution of some fractional order differential equations with the boundary conditions in Caputo sense. The main advantage of Caputo fractional derivative is the fractional derivative of constant is always zero as well as the basic derivative.

The Riemann-Liouvelli Fractional Derivative of order $\alpha$ of the function $f(x)$ is given as:

$$D^\alpha f(x) = \frac{1}{\Gamma(n-\alpha)} \frac{d^n}{dx^n} \int_a^x (x-y)^{n-\alpha-1} f(y) dy \quad \text{where } \alpha > 0.$$

The Caputo Fractional Derivative of order $\alpha$ of the function $f(x)$ is given as:

$$D^\alpha f(x) = \frac{1}{\Gamma(n-\alpha)} \int_a^x (x-y)^{n-\alpha-1} f^{(n)}(y) dy \quad \text{where } \alpha > 0.$$

Weight Functions: Throughout this research work we use weight functions as the modified Legendre polynomial of degree $n$ [11]:

$$p_n(x) = \left[\frac{1}{n!} \frac{d^n}{dx^n} (x^2 - x)^n - (-1)^n\right](x-1), \; p_n(0) = p_n(1) = 0, \; n \geq 1.$$

Similarly, the modified Bernoulli polynomial of degree $n$ as [11]: $B_0(x) = 1$ and

$$B_m(x) = \sum_{n=0}^{m} \frac{1}{n+1} \sum_{k=0}^{n} (-1)^k \binom{n}{k} (x+k)^m - \sum_{n=0}^{m} \frac{1}{n+1} \sum_{k=0}^{n} (-1)^k \binom{n}{k} (k)^m, \; m \geq 1.$$

**Weighted Residual Method**

The weighted residuals method is an approximation technique for solving boundary value problems that enrol trial functions satisfying the given boundary conditions and an integral formulation to minimize error, in normal sense, over the problem domain.

Given a fractional differential equation of the general form:

$$D^\alpha[y(x), x] = 0 \quad a < x < b, \tag{1}$$

subject to the homogeneous boundary conditions

$$y(a) = y(b) = 0. \tag{2}$$

The method of weighted residuals seeks an approximate solution of the form

$$\tilde{y}(x) = \sum_{i=1}^{n} a_i P_i(x) \tag{3}$$

where $\tilde{y}(x)$ denote the approximate solution can be expressed as the product of $a_i$ unknown, constant parameters to be determined and $P_i(x)$ are weight functions. The major requirement allocated on the trail functions treat as the permissible functions which are continuous over the domain and satisfy the specified boundary conditions. The residual function is also a function of unknown parameters $a_i$ and it can be expressed by

$$R(x) = D^\alpha[y(x), x] \neq 0. \tag{4}$$

The method of weighted residuals requires that unknown parameters $a_i$ be evaluated such that

$$\int_a^b P_j(x) R(x) dx = 0 \quad j = 1, 2, \ldots, n. \tag{5}$$

This equation can be solved for the $n$ values of $a_i$.

## 3. Formulation of nonlinear fractional order differential equations

In this section, we describe elaborately three weighted residual methods, namely, Galerkin method, Least-Square method and Collocation method, subsequently.

### (a) Formulation by Galerkin Method

To obtain the approximate solutions to the boundary value BVP in the Galerkin weighted residual (GWR) method using modified Legendre and Bernoulli polynomials as basis functions, we shall denote approximate trial solution by $\tilde{u}(x)$ by considering $u(x)$ denotes the exact solution to a boundary value problem. By supposing the nonlinear fractional order two-point boundary value problems with the boundary conditions,

$$\frac{d}{dx}\left(p(x)\frac{du}{dx}\right) + s(x)\frac{d^\alpha u}{dx^\alpha} + u^\beta(x) = f(x), \tag{6a}$$

$$u(a) = a_0,\ u(b) = b_0;\ \alpha \geq 1.5 \text{ and } \beta \geq 2. \tag{6b}$$

We assume the closeness solution of the differential equation as:

$$\tilde{u}(x) = \theta_0(x) + \sum_{j=1}^{n} a_j \theta_j(x). \tag{7}$$

Choose $\theta_0(x) = 0$ and $\theta_j(0) = \theta_j(1) = 0$ for each $j = 0, 1, \dots, n$.

Now the residual function is given by

$$\varepsilon(x) = \frac{d}{dx}\left(p(x)\frac{d\tilde{u}}{dx}\right) + s(x)\frac{d^\alpha \tilde{u}}{dx^\alpha} + \tilde{u}^\beta(x) - f(x). \tag{8}$$

The Galerkin weighted residual equations are then

$$\int_0^1 \varepsilon(x)\theta_i(x)dx = 0$$

or, $\int_0^1 \left[p(x)\frac{d^2\tilde{u}}{dx^2} + s(x)\frac{d^\alpha \tilde{u}}{dx^\alpha} + \tilde{u}^\beta(x)\right]\theta_i(x)dx - \int_0^1 f(x)\theta_i(x)dx = 0.$

or, equivalently

$$\int_0^1 \left[p(x)\frac{d^2\tilde{u}}{dx^2} + s(x)\frac{d^\alpha \tilde{u}}{dx^\alpha} + \tilde{u}^\beta(x)\right]\theta_i(x)dx = \int_0^1 f(x)\theta_i(x)dx. \tag{9}$$

Now

$$\int_0^1 p(x)\frac{d^2\tilde{u}}{dx^2}\theta_i(x)dx = -\int_0^1 p(x)\frac{d\tilde{u}}{dx}\frac{d\theta_i}{dx}dx \quad \text{since} \quad \theta_i(1) = \theta_i(0) = 0.$$

Equation (9) becomes

$$-\int_0^1 p(x)\frac{d\tilde{u}}{dx}\frac{d\theta_i}{dx}dx + \int_0^1 s(x)\frac{d^\alpha \tilde{u}}{dx^\alpha}\theta_i(x)dx + \int_0^1 \tilde{u}^\beta(x)\theta_i(x)dx = \int_0^1 f(x)\theta_i(x)dx. \tag{10}$$

Insert equation (2) into (10) and solving it, we get

$$\sum_{j=1}^{n}\left[-\int_0^1\left[p(x)\frac{d\theta_j}{dx}\frac{d\theta_i}{dx}+s(x)\frac{d^\alpha}{dx^\alpha}\theta_j(x)\theta_i(x)+(\theta_j(x))^\beta(x)\theta_i(x)\right]dx\right]=\int_0^1 f(x)\theta_i(x)dx.$$

which can be written as

$$\sum_{j=1}^{n}k_{i,j}a_j = F_i, \qquad (11)$$

where

$$k_{i,j} = -\int_0^1\left[p(x)\frac{d\theta_j}{dx}\frac{d\theta_i}{dx}+s(x)\frac{d^\alpha}{dx^\alpha}\theta_j(x)\theta_i(x)+(\theta_j(x))^\beta(x)\theta_i(x)\right]dx,$$

and $F_i = \int_0^1 f(x)\theta_i(x)dx$,

which is clearly the matrix form of a system of $n$ nonlinear equations.

Solving the system (11) yields the values of parameters and, upon substituting into equation (7) the approximate solution of the desired FBVP (6) is obtained.

### (b) Formulation by Least-Square Method

By considering a basis functions as modified Bernoulli and Legendre polynomials in the Least-Square Method, we obtain the approximate solutions to the boundary value problems.

If $v(x)$ is the exact solution to a boundary value problem, and then by denoting the approximate trial solution by $\tilde{v}(x)$, consider the nonlinear fractional order two-point BVP with the boundary conditions

$$p(x)\frac{dv}{dx}+s(x)\frac{d^\alpha v}{dx^\alpha}+v^\beta(x)=f(x), \qquad (12a)$$

$$v(a)=a_0,\ v(b)=b_0;\ \alpha \geq 1.5\ \text{and}\ \beta \geq 2. \qquad (12b)$$

We assume the closeness of solution of differential equation as:

$$\tilde{v}(x)=\varphi_0(x)+\sum_{j=1}^{n}a_j\varphi_j(x) \qquad (13)$$

Choose $\varphi_0(x)=0$ and $\varphi_j(0)=\varphi_j(1)=0$ for each $j=0,1,\ldots,n$.

Now the residual function is given by

$$R(x)=p(x)\frac{d\tilde{v}}{dx}+s(x)\frac{d^\alpha \tilde{v}}{dx^\alpha}+\tilde{v}^\beta(x)-f(x)=0. \qquad (14)$$

In this case, weight function is chosen as $W_j = \frac{\partial R(x)}{\partial a_j},\ j=1,2,\ldots,n.$

Now this choice of $W_j$ corresponds to minimize the mean square residual

$$WR=\frac{1}{2}\int R^2 dx = \text{minimum}.$$

The necessary condition for $WR$ to be minimum are given by

$$\frac{\partial WR}{\partial a_j} = \int_0^1 R(x) \frac{\partial R}{\partial a_j} dx = 0, \quad j = 1,2,\dots,n. \tag{15}$$

The equation (15) is clearly the matrix form of a system of $n$ nonlinear equations consisting of parameters $a_j$.

Solving the system (15), the values of parameters are determined and substitute into equation (13), the approximate solutions of the desired FBVP (12) are achieved.

(c) **Formulation by Collocation Method**

In this case, the approximate trial solution as $\widetilde{w}(x)$ where $w(x)$ is the unknown exact solution to a BVP, and consider nonlinear fractional order two-point BVPs with the boundary conditions:

$$p(x)\frac{dw}{dx} + s(x)\frac{d^\alpha w}{dx^\alpha} + w^\beta(x) = f(x), \tag{16a}$$

$$w(a) = a_0, \ w(b) = b_0; \ \alpha \geq 1.5 \text{ and } \beta \geq 2. \tag{16b}$$

We assume the closeness solution as:

$$\widetilde{w}(x) = \psi_0(x) + \sum_{j=1}^n a_j \psi_j(x). \tag{17}$$

Choose $\psi_0(x) = 0$ and $\psi_j(0) = \psi_j(1) = 0$ for each $j = 0, 1, \dots, n$, $a_j$ are the unknown parameters and $\psi_j(x)$ are the basis functions. We choose modified Legendre and Bernoulli polynomials as basis functions. In another case, $\psi_0(x)$ is defined to satisfy the nonhomogeneous boundary condition so that other basis functions satisfy the homogeneous boundary conditions.

Now the residual function is given by

$$R(x) = p(x)\frac{d\widetilde{w}}{dx} + s(x)\frac{d^\alpha \widetilde{w}}{dx^\alpha} + \widetilde{w}^\beta(x) - f(x). \tag{18}$$

Insert the equation (17) into the equation (16a), we obtain the equation of the form:

$$R(x) = p(x)\frac{d}{dx}\left(\psi_0(x) + \sum_{j=1}^n a_j \psi_j(x)\right) + s(x)\frac{d^\alpha}{dx^\alpha}\left(\psi_0(x) + \sum_{j=1}^n a_j \psi_j(x)\right) +$$

$$\left(\psi_0(x) + \sum_{j=1}^n a_j \psi_j(x)\right)^\beta (x) - f(x).$$

In Collocation method, we evaluate the residual function at some grid points $x_j$ and setting the residual function as $R(x_j) = 0$ by the arrangement of fractional differential equation with the boundary conditions.

We assume that the boundary conditions on $[a,b]$ such as $w(a) = a_0, \ w(b) = b_0$. If we choose the $n$ parameters and the boundary point starts from $a$ then the grid points are described as:

$$x_j = \frac{a+j}{n+1}, \quad \text{where } j = 1,2,3,\ldots,n.$$

Setting $R(x_j) = 0$, we obtain the system in unknown parameters $a_j$. Putting the values of the parameters into equation (17), we get the approximate solution of nonlinear fractional order BVP (16).

## 4. Test Problems

In this section, we test the proposed formulations by considering numerical examples which are available in the literature and investigated by several methods. Further, the efficiency and reliability of the proposed method is applied for these problems by computing the error $L_2$ and the maximum absolute error $L_\infty$ which are given as follows:

$$L_2 = \int_0^1 (u(x) - \tilde{u}(x))dx \text{ and } L_\infty = max|u(x) - \tilde{u}(x)|,$$

where $u(x)$ and $\tilde{u}(x)$ are the exact and approximate solutions, respectively.

**Problem 1:** Consider the non-linear fractional BVP [18]:

$$D^{1.5}u(x) - u^3(x) = f(x) \text{ with boundary conditions } u(0) = -1, \ u(1) = 0. \quad (19)$$

where $f(x) = \frac{\Gamma(2.9)}{\Gamma(1.4)}x^{0.4} - (x^{1.9} - 1)^3$.

The exact solution of this problem is $u(x) = x^{1.9} - 1$. The approximate solution is derived with respect to the unknown coefficients $a_j$ from the equation (2), (8) and (12). Three weighted residual methods Galerkin, Least Square and Collocation give the approximate solution $\tilde{u}(x), \tilde{v}(x)$ and $\tilde{w}(x)$, respectively, of the given problem using the modified Legendre polynomial of degree $(n = 3)$ as basis functions, we have:

$$\tilde{u}_{GM}(x) = -1 + 0.01506978x + 1.1353053x^2 - 0.2236533x^3 + 0.0732781x^4,$$

$$\tilde{v}_{LSM}(x) = -1 + 0.0128251x + 1.1408130x^2 - 0.2227155x^3 + 0.0690774x^4,$$

$$\tilde{w}_{CM}(x) = -1 + 0.0193695x + 1.1142138x^2 - 0.1962006x^3 + 0.0626172x^4.$$

Similarly, when we use the modified Bernoulli polynomial as basis functions of degree $(n = 3)$ we get another approximate solution as given below:

$$\tilde{u}_{GM}(x) = -1 + 0.0151999x + 1.1351020x^2 - 0.22341003x^3 + 0.0732156x^4,$$

$$\tilde{v}_{LSM}(x) = -1 + 0.0128251x + 1.1408130x^2 - 0.2227155x^3 + 0.0690774x^4,$$

$$\tilde{w}_{CM}(x) = -1 + 0.0193800x + 1.1141910x^2 - 0.1961811x^3 + 0.0626100x^4.$$

**Table 1.1:** Absolute errors of the BVP in Equation (19)

| $x$ | Exact Solutions | Absolute errors obtained using modified Legendre polynomials of degree ($n = 3$) | | | Absolute errors obtained using modified Bernoulli polynomials of degree ($n = 3$) | | |
|---|---|---|---|---|---|---|---|
| | | GWR | Least-Square | Collocation | GWR | Least-Square | Collocation |
| 0   | $-1$    | 0 | 0 | 0 | 0 | 0 | 0 |
| 0.1 | $-0.987$ | $5.44 \times 10^{-5}$ | $2.12 \times 10^{-3}$ | $2.99 \times 10^{-4}$ | $6.46 \times 10^{-5}$ | $1.14 \times 10^{-4}$ | $3.00 \times 10^{-4}$ |
| 0.2 | $-0.953$ | $2.30 \times 10^{-4}$ | $2.19 \times 10^{-3}$ | $1.17 \times 10^{-5}$ | $2.14 \times 10^{-4}$ | $4.58 \times 10^{-4}$ | $1.04 \times 10^{-5}$ |
| 0.3 | $-0.898$ | $2.61 \times 10^{-4}$ | $1.58 \times 10^{-3}$ | $2.15 \times 10^{-4}$ | $2.43 \times 10^{-4}$ | $4.48 \times 10^{-4}$ | $2.13 \times 10^{-4}$ |
| 0.4 | $-0.824$ | $1.14 \times 10^{-4}$ | $8.16 \times 10^{-4}$ | $2.85 \times 10^{-4}$ | $9.70 \times 10^{-5}$ | $1.78 \times 10^{-4}$ | $2.83 \times 10^{-4}$ |
| 0.5 | $-0.732$ | $4.10 \times 10^{-5}$ | $1.50 \times 10^{-4}$ | $3.16 \times 10^{-4}$ | $5.66 \times 10^{-5}$ | $1.50 \times 10^{-4}$ | $3.15 \times 10^{-4}$ |
| 0.6 | $-0.621$ | $7.20 \times 10^{-5}$ | $2.97 \times 10^{-4}$ | $3.93 \times 10^{-4}$ | $8.48 \times 10^{-5}$ | $3.66 \times 10^{-4}$ | $3.91 \times 10^{-4}$ |
| 0.7 | $-0.492$ | $6.30 \times 10^{-5}$ | $4.92 \times 10^{-4}$ | $5.31 \times 10^{-4}$ | $5.33 \times 10^{-5}$ | $3.77 \times 10^{-4}$ | $5.30 \times 10^{-4}$ |
| 0.8 | $-0.345$ | $2.86 \times 10^{-4}$ | $4.60 \times 10^{-4}$ | $6.55 \times 10^{-4}$ | $2.79 \times 10^{-4}$ | $2.02 \times 10^{-4}$ | $6.55 \times 10^{-4}$ |
| 0.9 | $-0.181$ | $3.84 \times 10^{-4}$ | $2.64 \times 10^{-4}$ | $5.80 \times 10^{-4}$ | $3.81 \times 10^{-4}$ | $1.61 \times 10^{-5}$ | $5.80 \times 10^{-5}$ |
| 1   | 0       | 0 | 0 | 0 | 0 | 0 | 0 |

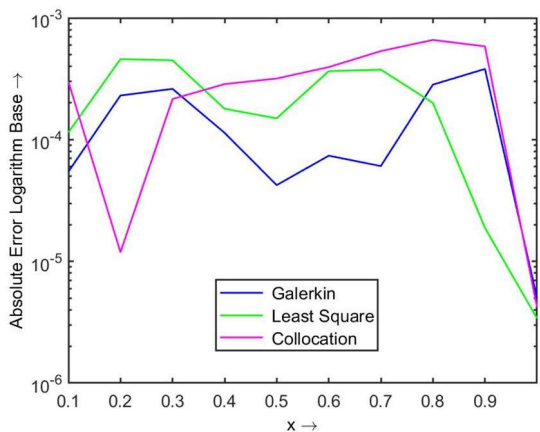
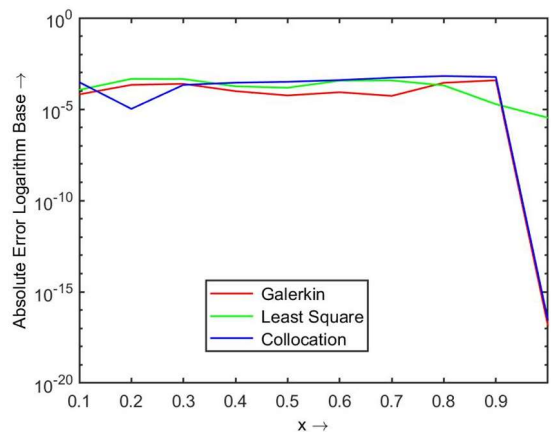

**Fig1.1:** Absolute errors using modified Legendre polynomial     **Fig1.2:** Absolute error using modified Bernoulli polynomial

To find the accuracy of the approximate solution, we compare the result with the exact solution and the absolute error can be defined by $|u_{Exact}(x) - \tilde{u}_{Appx}(x)|$. Then the obtained absolute error using two different polynomials by three residual methods are shown in Table 1.1, and the comparison of their results are displayed in Figs. 1.1 and 1.2.

**Table 1.2:** $L_\infty$ and $L_2$ error of the Problem1 in Eqn. (19)

| Errors | Errors obtained using modified Legendre polynomials of degree ($n = 3$) | | | Errors obtained using modified Bernoulli polynomials degree ($n = 3$) | | | Reference [18] |
|---|---|---|---|---|---|---|---|
| | GWR | Least Square | Collocation | GWR | Least Square | Collocation | New hybrid method ($N = 10$) |
| $L_\infty$ | $3.84 \times 10^{-4}$ | $2.19 \times 10^{-3}$ | $6.55 \times 10^{-4}$ | $3.81 \times 10^{-4}$ | $4.58 \times 10^{-4}$ | $6.55 \times 10^{-4}$ | $5.66 \times 10^{-5}$ |
| $L_2$ | $4.02 \times 10^{-8}$ | $8.08 \times 10^{-8}$ | $1.56 \times 10^{-7}$ | $3.80 \times 10^{-8}$ | $7.49 \times 10^{-8}$ | $1.55 \times 10^{-7}$ | $8.92 \times 10^{-6}$ |

From Table 1.1 and Table 1.2 we observe that the solutions converge to the exact solutions monotonically, and good agreement with the exact solutions even using lower order polynomials as weight functions.

**Problem 2:** Suppose the non-linear fractional BVP [16]

$$D^{1.5}u(x) = -2(u'(x))^2 - 8u(x), \quad 0 \le x \le 1,$$
$$u(0) = 0, \quad u'(1) = -1. \tag{20}$$

The exact solution of this problem is $u(x) = x - x^2$. The approximate solution $\tilde{u}(x), \tilde{v}(x)$ and $\tilde{w}(x)$ of the given problem using the modified Legendre polynomial of degree ($n = 3$) as basis functions:

$\tilde{u}_{GM}(x) = x - x^2 + 1.87875 \times 10^{-12}x^3 - 8.88285 \times 10^{-13}x^4,$

$\tilde{v}_{LSM}(x) = x - x^2 + 8.47192 \times 10^{-12}x^3 - 5.60912 \times 10^{-12}x^4,$

$\tilde{w}_{CM}(x) = 0.99999x - 0.99999x^2 - 1.74997 \times 10^{-15}x^3 + 6.65225 \times 10^{-16}x^4.$

Similarly, when we use the modified Bernoulli polynomial of degree ($n = 3$) as a basis function we get another approximate solution as given below:

$\tilde{u}_{GM}(x) = 0.99999x - 0.99999x^2 - 1.99878 \times 10^{-1} \, x^3 + 8.88789 \times 10^{-1} \, x^4,$

$\tilde{v}_{LSM}(x) = 0.99999x - 0.99999x^2 - 2.87383 \times 10^{-1} \, x^3 + 1.86700 \times 10^{-12}x^4,$

$\tilde{w}_{CM}(x) = 0.99999x - 0.99999x^2 - 1.92650 \times 10^{-1} \, x^3 + 6.46201 \times 10^{-16}x^4.$

The obtained absolute errors using two different polynomials by three residual methods are shown in the Table 2.1, and the comparison are displayed in the Figs. 2.1 and 2.2.

**Table 2.1:** Absolute error for problem in Eqn. (20)

| $x$ | Exact Solutions | Absolute errors obtained using Modified Legendre polynomials ($n = 3$) | | | Absolute errors obtained using Modified Bernoulli polynomials ($n = 3$) | | |
|---|---|---|---|---|---|---|---|
| | | GWRM | Least Square | Collocation | GWRM | Least Square | Collocation |
| 0 | 0 | 0 | 0 | 0 | 0 | 0 | 0 |
| 0.1 | 0.09 | $1.70 \times 10^{-14}$ | $1.30 \times 10^{-14}$ | $4.16 \times 10^{-17}$ | $2.21 \times 10^{-14}$ | $4.14 \times 10^{-15}$ | $5.55 \times 10^{-17}$ |
| 0.2 | 0.16 | $1.87 \times 10^{-14}$ | $4.20 \times 10^{-14}$ | $2.77 \times 10^{-17}$ | $2.58 \times 10^{-14}$ | $1.45 \times 10^{-14}$ | $5.55 \times 10^{-17}$ |
| 0.3 | 0.21 | $1.31 \times 10^{-14}$ | $5.64 \times 10^{-14}$ | $2.77 \times 10^{-17}$ | $1.99 \times 10^{-14}$ | $2.06 \times 10^{-14}$ | $8.32 \times 10^{-17}$ |
| 0.4 | 0.24 | $6.13 \times 10^{-15}$ | $3.89 \times 10^{-14}$ | $2.77 \times 10^{-17}$ | $1.11 \times 10^{-14}$ | $1.63 \times 10^{-14}$ | $8.32 \times 10^{-17}$ |
| 0.5 | 0.25 | $1.61 \times 10^{-15}$ | $1.41 \times 10^{-14}$ | $2.77 \times 10^{-17}$ | $3.83 \times 10^{-15}$ | $2.22 \times 10^{-16}$ | $5.55 \times 10^{-17}$ |
| 0.6 | 0.24 | $1.22 \times 10^{-15}$ | $9.29 \times 10^{-14}$ | 0 | $4.99 \times 10^{-16}$ | $2.49 \times 10^{-14}$ | $2.77 \times 10^{-17}$ |
| 0.7 | 0.20 | $4.52 \times 10^{-15}$ | $1.74 \times 10^{-13}$ | 0 | $1.38 \times 10^{-15}$ | $5.16 \times 10^{-14}$ | $2.77 \times 10^{-17}$ |
| 0.8 | 0.15 | $8.93 \times 10^{-15}$ | $2.21 \times 10^{-13}$ | 0 | $4.60 \times 10^{-15}$ | $6.80 \times 10^{-14}$ | 0 |
| 0.9 | 0.08 | $9.70 \times 10^{-15}$ | $1.84 \times 10^{-13}$ | 0 | $6.20 \times 10^{-15}$ | $5.77 \times 10^{-14}$ | 0 |
| 1 | 0 | 0 | 0 | 0 | 0 | 0 | 0 |

**Table 2.2:** $L_\infty$ and $L_2$ errors of the problem in Eqn. (20)

| Errors | Errors obtained using modified Legendre polynomials of degree ($n = 3$) | | | Errors obtained using modified Bernoulli polynomials of degree ($n = 3$) | | |
|---|---|---|---|---|---|---|
| | Galerkin | Least-Square | Collocation | Galerkin | Least-Square | Collocation |
| $L_\infty$ | $1.87 \times 10^{-14}$ | $2.21 \times 10^{-13}$ | $4.16 \times 10^{-17}$ | $2.58 \times 10^{-14}$ | $6.80 \times 10^{-14}$ | $8.32 \times 10^{-17}$ |
| $L_2$ | $1.05 \times 10^{-28}$ | $1.29 \times 10^{-26}$ | $8.30 \times 10^{-34}$ | $1.76 \times 10^{-28}$ | $1.22 \times 10^{-27}$ | $1.76 \times 10^{-33}$ |

Table 2.1 and Table 2.2 reflect that the solutions converge fast to the exact solutions, and a very good agreement with the exact solutions even using lower degree polynomials as weight functions.

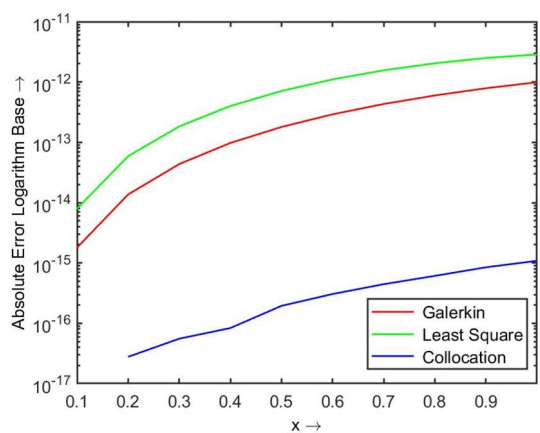 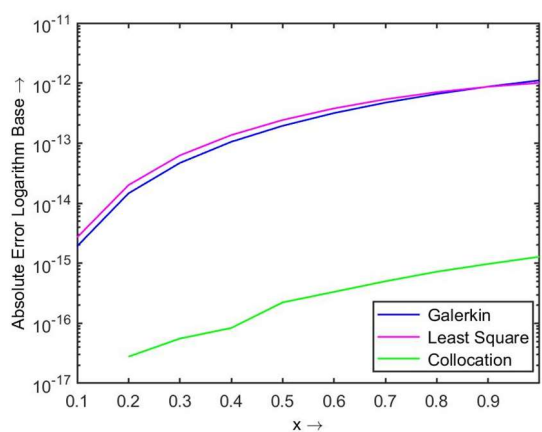

**Fig2.1:** Absolute error using modified Legendre polynomial    **Fig2.2:** Absolute error using modified Bernoulli polynomial

**Example 3:** Consider the non-linear fractional BVP [16]:

$$D^\alpha u = u^2(x) + 90.27 x^{1.5} - 2 \times 27.08 x^{0.5} - 2 \times \frac{0.56}{x^{1.5}} - x^{10}$$
$$+ 4x^9 - 4x^8 - 4x^7 + 8x^6 - 4x^4, \quad (21)$$
$$u(0) = 0, \quad u'(0) = 0, \quad u(1) = 1, \quad u'(1) = 1,$$

with $\alpha = 3.5$, $0 \leq x \leq 1$ and $3 \leq \alpha \leq 4$.

The exact solution of this problem is $u(x) = x^5 - 2x^4 + 2x^2$. The approximate solution $\tilde{u}(x), \tilde{v}(x)$ and $\tilde{w}(x)$ of the given problem using the modified Legendre polynomial of degree ($n = 5$) as basis functions are as follows:

$$\tilde{u}_{GM}(x) = 6.55 \times 10^{-1}\ x + 1.99 x^2 + 1.98 \times 10^{-1}\ x^3 - 2x^4 + x^5 - 1.96 \times 10^{-1}\ x^6,$$
$$\tilde{v}_{LSM}(x) = 1.11 \times 10^{-16} x + 1.99 x^2 + 1.88 \times 10^{-11} x^3 - 2x^4 + x^5 - 1.96 \times 10^{-1}\ x^6,$$
$$\tilde{w}_{CM}(x) = -4.55 \times 10^{-15} x + 1.99 x^2 + 3.41 \times 10^{-14} x^3 - 2x^4 + x^5 - 3.81 \times 10^{-15} x^6.$$

while using the modified Bernoulli polynomial as basis functions we get the approximations as given below:

$$\tilde{u}_{GM}(x) = -2.38 \times 10^{-14} x + 1.99 x^2 + 6.19 \times 10^{-12} x^3 - 2x^4 + x^5 - 1.95 \times 10^{-12} x^6,$$
$$\tilde{v}_{LSM}(x) = 6.66 \times 10^{-13} x + 1.99 x^2 + 7.33 \times 10^{-12} x^3 - 2x^4 + x^5 - 1.95 \times 10^{-12} x^6,$$
$$\tilde{w}_{CM}(x) = 3.38 \times 10^{-15} x + 2x^2 - 2.15 \times 10^{-14} x^3 - 1.99 x^4 + x^5 - 1.93 \times 10^{-15} x^6.$$

The graphical representation of the two solutions is delineated in the Figs. 3.1 and 3.2, which shows that the approximate solution is in sensible agreement with the exact solution. The differences in the exact and approximate solutions are scarcely perceivable.

From Figures 3.1 and 3.2, and Table 3.1 we may notice that the solutions converge fast to the exact solutions, and a very good agreement with the exact solutions on using polynomials of degree ($n = 5$) as weight functions.

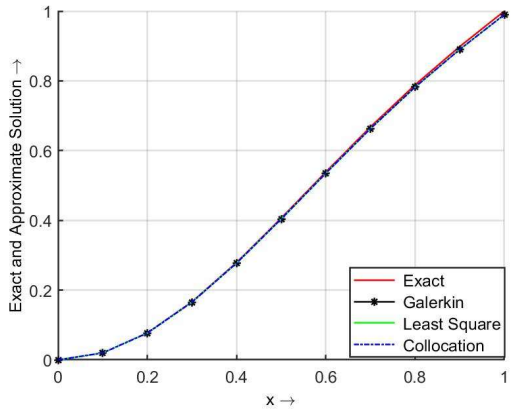 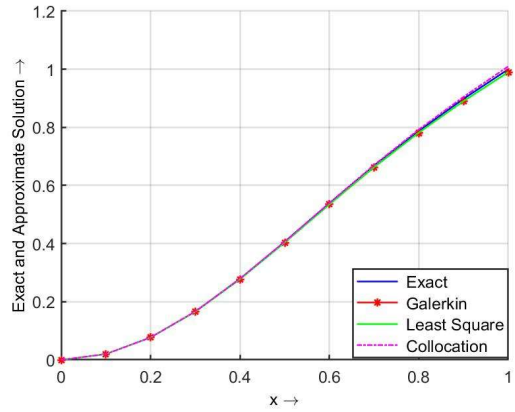

**Fig 3.1**: Exact and approximate solutions of Eqn. (21) using modified Legendre polynomials

**Fig 3.2**: Exact and approximate solutions of Eqn. (21) using modified Bernoulli polynomials

**Table 3.1:** Absolute errors of problem 3 in Eqn. (21)

| $x$ | Exact Solutions | Absolute errors obtained by existing method using modified Legendre polynomials ($n = 5$ and $\alpha = 3.5$) | | | Absolute errors obtained by existing method using modified Bernoulli polynomials ($n = 5$ and $\alpha = 3.5$) | | |
|---|---|---|---|---|---|---|---|
| | | Galerkin | Least-Square | Collocation | Galerkin | Least-Square | Collocation |
| 0 | 0 | 0 | 0 | 0 | 0 | 0 | 0 |
| 0.1 | 0.019 | $4.26 \times 10^{-14}$ | $3.71 \times 10^{-15}$ | $5.89 \times 10^{-16}$ | $1.54 \times 10^{-14}$ | $3.78 \times 10^{-14}$ | $4.26 \times 10^{-16}$ |
| 0.2 | 0.077 | $9.27 \times 10^{-14}$ | $3.95 \times 10^{-14}$ | $1.26 \times 10^{-15}$ | $4.16 \times 10^{-14}$ | $3.81 \times 10^{-14}$ | $9.57 \times 10^{-16}$ |
| 0.3 | 0.166 | $1.48 \times 10^{-13}$ | $9.83 \times 10^{-14}$ | $1.88 \times 10^{-15}$ | $6.53 \times 10^{-14}$ | $2.01 \times 10^{-14}$ | $1.47 \times 10^{-15}$ |
| 0.4 | 0.279 | $1.80 \times 10^{-13}$ | $1.44 \times 10^{-13}$ | $2.55 \times 10^{-15}$ | $8.22 \times 10^{-14}$ | $5.77 \times 10^{-15}$ | $1.88 \times 10^{-15}$ |
| 0.5 | 0.406 | $1.67 \times 10^{-13}$ | $1.48 \times 10^{-13}$ | $2.88 \times 10^{-15}$ | $9.19 \times 10^{-14}$ | $3.34 \times 10^{-14}$ | $2.10 \times 10^{-15}$ |
| 0.6 | 0.538 | $1.07 \times 10^{-13}$ | $1.05 \times 10^{-15}$ | $2.88 \times 10^{-15}$ | $9.48 \times 10^{-14}$ | $5.82 \times 10^{-14}$ | $2.22 \times 10^{-15}$ |
| 0.7 | 0.667 | $2.79 \times 10^{-14}$ | $3.54 \times 10^{-14}$ | $2.99 \times 10^{-15}$ | $8.99 \times 10^{-14}$ | $7.42 \times 10^{-14}$ | $2.10 \times 10^{-15}$ |
| 0.8 | 0.788 | $2.94 \times 10^{-14}$ | $1.94 \times 10^{-14}$ | $2.22 \times 10^{-15}$ | $7.46 \times 10^{-14}$ | $7.46 \times 10^{-14}$ | $1.66 \times 10^{-15}$ |
| 0.9 | 0.898 | $3.53 \times 10^{-14}$ | $2.75 \times 10^{-14}$ | $1.44 \times 10^{-15}$ | $4.54 \times 10^{-14}$ | $5.17 \times 10^{-14}$ | $8.88 \times 10^{-16}$ |
| 1 | 1 | 0 | 0 | 0 | 0 | 0 | 0 |
| $L_\infty$ | | $1.80 \times 10^{-13}$ | $1.48 \times 10^{-13}$ | $2.99 \times 10^{-15}$ | $9.48 \times 10^{-14}$ | $7.46 \times 10^{-14}$ | $2.22 \times 10^{-15}$ |

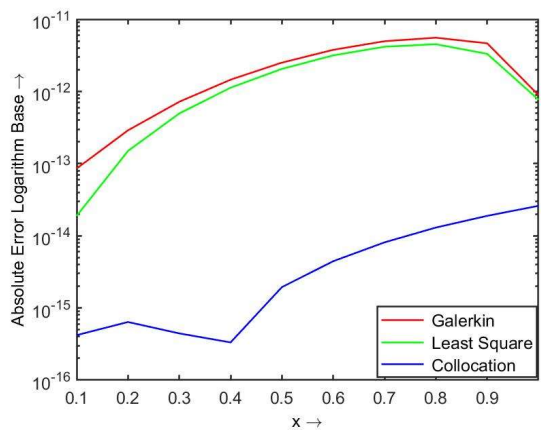 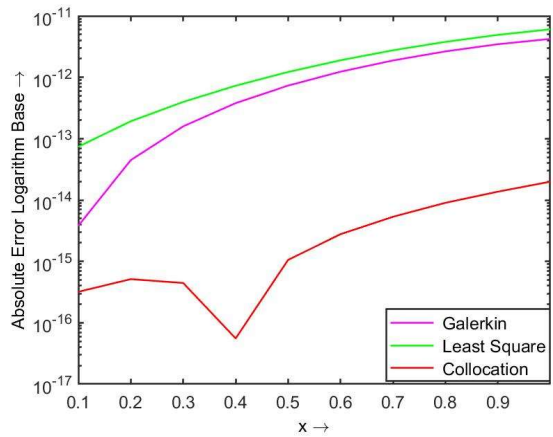

**Fig3.3:** Absolute error using modified Legendre polynomial

**Fig3.4:** Absolute error using modified Bernoulli polynomial

**Problem 4**: Consider the following FDE of the non-linear form [17]:

$$D^2 u(x) + \Gamma\left(\frac{4}{5}\right) x^{\frac{6}{5}} D^{\frac{6}{5}} u(x) + \frac{11}{9} \Gamma\left(\frac{5}{6}\right) x^{\frac{1}{6}} D^{\frac{1}{6}} u(x) - (Du(x))^2 = 2 + \frac{1}{10} x^2, \quad (22)$$
$$x(0) = 1, x(1) = 2.$$

The exact solution of Eqn. (22) is: $u(x) = 1 + x^2$.

The approximate solutions $\tilde{u}(x), \tilde{v}(x)$ and $\tilde{w}(x)$ of the given problem using the modified Legendre polynomial of degree ($n = 3$) as basis functions are:

$$\tilde{u}_{GM}(x) = 1 - 0.00133x + 1.00029x^2 - 0.00058x^3 + 0.00162x^4,$$
$$\tilde{v}_{LSM}(x) = 1 - 0.00128x + 1.00023x^2 - 0.00064x^3 + 0.00170x^4,$$
$$\tilde{w}_{CM}(x) = 1 - 0.00133x + 1.00022x^2 - 0.00033x^3 + 0.00144x^4,$$

Similarly, when we use the modified Bernoulli polynomial of degree ($n = 3$) as basis functions we get approximate solutions as given below:

$$\tilde{u}_{GM}(x) = 1 - 0.00133x + 1.00029x^2 - 0.00059x^3 + 0.00162x^4,$$
$$\tilde{v}_{LSM}(x) = 1 - 0.00128x + 1.00022x^2 - 0.00064x^3 + 0.00170x^4,$$
$$\tilde{w}_{CM}(x) = 1 - 0.00139x + 1.00020x^2 - 0.00031x^3 + 0.00144 x^4.$$

**Table 4.1:** Absolute errors of the problem 4 in Eqn. (22)

| $x$ | Absolute errors obtained using modified Legendre polynomials of degree ($n = 3$) | | | Absolute errors obtained using modified Bernoulli polynomials of degree ($n = 3$) | | | Reference |
|---|---|---|---|---|---|---|---|
| | GWRM | Least Square | Collocation | GWRM | Least Square | Collocation | Legendre Wavelet Operation Matrix Method |
| 0 | 0 | 0 | 0 | 0 | 0 | 0 | $2.00 \times 10^{-9}$ |
| 0.1 | $1.30 \times 10^{-4}$ | $1.27 \times 10^{-4}$ | $1.31 \times 10^{-4}$ | $1.27 \times 10^{-4}$ | $1.26 \times 10^{-4}$ | $1.30 \times 10^{-4}$ | $4.26 \times 10^{-4}$ |
| 0.2 | $2.57 \times 10^{-4}$ | $2.51 \times 10^{-4}$ | $2.59 \times 10^{-4}$ | $2.50 \times 10^{-4}$ | $2.49 \times 10^{-4}$ | $2.56 \times 10^{-4}$ | $7.58 \times 10^{-4}$ |
| 0.3 | $3.76 \times 10^{-4}$ | $3.69 \times 10^{-4}$ | $3.78 \times 10^{-4}$ | $3.66 \times 10^{-4}$ | $3.67 \times 10^{-4}$ | $3.75 \times 10^{-4}$ | $9.95 \times 10^{-4}$ |
| 0.4 | $4.82 \times 10^{-4}$ | $4.76 \times 10^{-4}$ | $4.84 \times 10^{-4}$ | $4.69 \times 10^{-4}$ | $4.74 \times 10^{-4}$ | $4.79 \times 10^{-4}$ | $1.13 \times 10^{-5}$ |
| 0.5 | $5.65 \times 10^{-4}$ | $5.61 \times 10^{-4}$ | $5.64 \times 10^{-4}$ | $5.50 \times 10^{-4}$ | $5.58 \times 10^{-4}$ | $5.60 \times 10^{-4}$ | $1.18 \times 10^{-5}$ |
| 0.6 | $6.10 \times 10^{-4}$ | $6.08 \times 10^{-4}$ | $6.07 \times 10^{-4}$ | $5.94 \times 10^{-4}$ | $6.06 \times 10^{-4}$ | $6.02 \times 10^{-4}$ | $1.13 \times 10^{-5}$ |
| 0.7 | $6.00 \times 10^{-4}$ | $6.01 \times 10^{-4}$ | $5.94 \times 10^{-4}$ | $5.84 \times 10^{-4}$ | $5.99 \times 10^{-4}$ | $5.90 \times 10^{-4}$ | $9.95 \times 10^{-4}$ |
| 0.8 | $5.13 \times 10^{-4}$ | $5.16 \times 10^{-4}$ | $5.05 \times 10^{-4}$ | $4.99 \times 10^{-4}$ | $5.14 \times 10^{-4}$ | $5.02 \times 10^{-4}$ | $7.58 \times 10^{-4}$ |
| 0.9 | $3.22 \times 10^{-4}$ | $3.25 \times 10^{-4}$ | $3.16 \times 10^{-4}$ | $3.14 \times 10^{-4}$ | $3.25 \times 10^{-4}$ | $3.14 \times 10^{-4}$ | $4.26 \times 10^{-4}$ |
| 1 | 0 | 0 | 0 | 0 | 0 | 0 | 0 |
| $L_\infty$ | $6.10 \times 10^{-4}$ | $6.08 \times 10^{-4}$ | $6.07 \times 10^{-4}$ | $5.94 \times 10^{-4}$ | $6.06 \times 10^{-4}$ | $6.02 \times 10^{-4}$ | $9.95 \times 10^{-4}$ |

From Table 4.1, and Figures 4.1 and 4.2, we may observe that the solutions converge monotonically to the exact solutions, and the accuracy agreement is considerable only using lower degree polynomials as weight functions. Finally, we may note that if we increase rapidly the degree of the polynomials then it may fail the monotonicity of the convergence, computational costs are high and the results may not stable.

**Conclusions**

In this research work, we have exploited the weighted residual methods, namely, Galerkin method, Least-Square method and Collocation method to find the approximate solutions to the fractional order non-linear boundary value problems rigorously. The computed results show that all the proposed three methods are effective, reliable and converge monotonically to the exact solutions. Finally, we conclude that these methods may be applied to find the approximate solutions to any kind of fractional order differential equations.

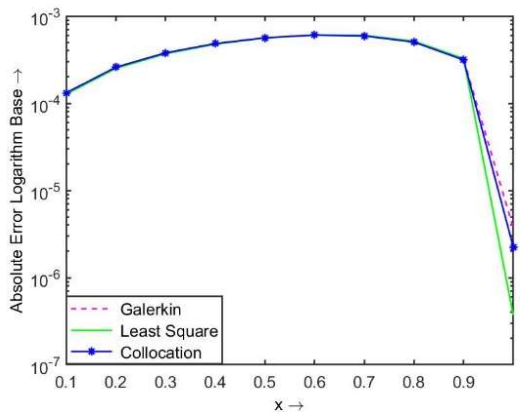 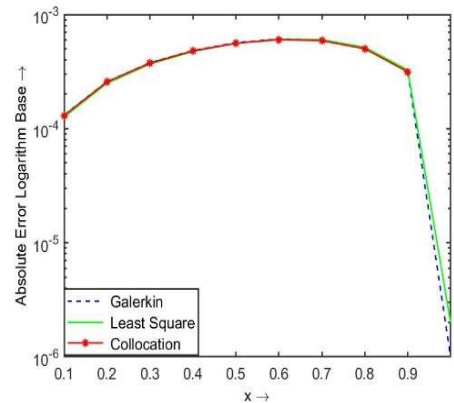

**Fig4.1:** Absolute error using modified Legendre polynomial    **Fig4.2:** Absolute error using modified Bernoulli polynomial